\documentclass[12pt,twoside]{article}

\usepackage{amssymb}
\usepackage{amsmath}
\usepackage{bbm}
\usepackage{mathrsfs}
\usepackage{float}
\usepackage{xypic}

\makeatletter
\def\hsmash{\relax 
  \ifmmode\def\next{\mathpalette\mathhsm@sh}\else\let\next\makehsm@sh
  \fi\next}
\def\makehsm@sh#1{\setbox\z@\hbox{#1}\finhsm@sh}
\def\mathhsm@sh#1#2{\setbox\z@\hbox{$\m@th#1{#2}$}\finhsm@sh}
\def\finhsm@sh{\wd\z@\z@ \box\z@}
\makeatother

\makeatletter
\gdef\th@mychange{\normalfont\slshape
   \def\@begintheorem##1##2{\item
        [\hskip\labelsep \theorem@headerfont ##2. ##1  \,--\!--\!--\!--  ]}%
 \def\@opargbegintheorem##1##2##3{%
   \item[\hskip\labelsep \theorem@headerfont ##2. ##1\ {\upshape(}##3{\upshape)}. \,-----  ]}}
\makeatother

\RequirePackage{theorem}
\theoremstyle{mychange}

{\theorembodyfont{\rmfamily}\newtheorem{ttt}{}[subsection]}
{\theorembodyfont{\rmfamily}}
{\theorembodyfont{\rmfamily}\newtheorem{rems}[ttt]{Remarks.}}
{\theorembodyfont{\rmfamily}}
{\theorembodyfont{\rmfamily}}
{\theorembodyfont{\rmfamily}}
{\theorembodyfont{\rmfamily}\newtheorem{summa}[ttt]{Summary.}}

{\theorembodyfont{\itshape}\newtheorem{fac}[ttt]{Fact.}}
{\theorembodyfont{\itshape}\newtheorem{lem}[ttt]{Lemma.}}
{\theorembodyfont{\itshape}\newtheorem{sub}[ttt]{Sublemma.}}
{\theorembodyfont{\itshape}\newtheorem{prop}[ttt]{Proposition.}}
{\theorembodyfont{\itshape}}
{\theorembodyfont{\itshape}\newtheorem{theorem}[ttt]{Theorem.}}
{\theorembodyfont{\itshape}}

{\theorembodyfont{\itshape}\newtheorem{theoremo}[ttt]{Theorem}}
{\theorembodyfont{\itshape}}
{\theorembodyfont{\itshape}}

{\theorembodyfont{\rmfamily}\newtheorem{ttts}{}[section]}
{\theorembodyfont{\rmfamily}\newtheorem{remes}[ttts]{Remark.}}
{\theorembodyfont{\rmfamily}\newtheorem{remss}[ttts]{Remarks.}}
{\theorembodyfont{\rmfamily}\newtheorem{defis}[ttts]{Definition.}}
{\theorembodyfont{\rmfamily}\newtheorem{defiss}[ttts]{Definitions.}}
{\theorembodyfont{\rmfamily}\newtheorem{exes}[ttts]{Example.}}
{\theorembodyfont{\rmfamily}}

{\theorembodyfont{\itshape}\newtheorem{facs}[ttts]{Fact.}}
{\theorembodyfont{\itshape}\newtheorem{lems}[ttts]{Lemma.}}
{\theorembodyfont{\itshape}}
{\theorembodyfont{\itshape}\newtheorem{coros}[ttts]{Corollary.}}

{\theorembodyfont{\rmfamily}\newtheorem{exess}[ttts]{Example}}
{\theorembodyfont{\rmfamily}\newtheorem{algos}[ttts]{Algorithm}}

{\theorembodyfont{\itshape}\newtheorem{facos}[ttts]{Fact}}
{\theorembodyfont{\itshape}}
{\theorembodyfont{\itshape}\newtheorem{lemos}[ttts]{Lemma}}

{\theorembodyfont{\rmfamily}\newtheorem{pls}[ttts]{Plan of the article.}}
{\theorembodyfont{\rmfamily}\newtheorem{dos}[ttts]{Description of the sample.}}
{\theorembodyfont{\rmfamily}\newtheorem{tres}[ttts]{The results.}}
{\theorembodyfont{\rmfamily}}
{\theorembodyfont{\rmfamily}\newtheorem{suma}[ttts]{Summary.}}

\newcounter{III}
\newenvironment{III}{\begin{list}{\rm \Roman{III}) }{\usecounter{III} \leftmargin=0.0pt \labelsep=0.0pt \listparindent=0.0pt \labelwidth=0.0pt \parsep=\smallskipamount \itemsep=0.0pt \topsep=0.0pt \partopsep=\smallskipamount}}{\end{list}}

\newcounter{iii}
\newenvironment{iii}{\begin{list}{\rm \roman{iii}) }{\usecounter{iii} \leftmargin=0.0pt \labelsep=0.0pt \listparindent=0.0pt \labelwidth=0.0pt \parsep=\smallskipamount \itemsep=0.0pt \topsep=0.0pt \partopsep=\smallskipamount}}{\end{list}}

\newcounter{abc}
\newenvironment{abc}{\begin{list}{\rm \alph{abc}) }%
{\usecounter{abc} \leftmargin=0.0pt \labelsep=0.0pt %
\listparindent=0.0pt \labelwidth=0.0pt \parsep=\smallskipamount %
\itemsep=0.0pt \topsep=0.0pt \partopsep=\smallskipamount}}{\end{list}}

\newcommand{\Pb}{{\mathop{\text{\bf P}}}}
\newcommand{\Ab}{{\mathop{\text{\bf A}}}}

\newcommand{\et}{\text{\rm {\'e}t}}
\newcommand{\ab}{\text{\rm ab}}
\newcommand{\geo}{\text{\rm geo}}
\newcommand{\tors}{\text{\rm tors}}
\newcommand{\Spec}{\mathop{\text{\rm Spec}}}

\newcommand{\Frob}{\text{\rm Frob}}
\newcommand{\Br}{\text{\rm Br}}

\newcommand{\Pic}{\mathop{\text{\rm Pic}}}

\newcommand{\disc}{\mathop{\text{\rm disc}}}
\renewcommand{\div}{\mathop{\text{\rm div}}}
\newcommand{\MW}{\mathop{\text{\rm MW}}}

\newcommand{\Z}{\mathop{\text{\rm Z}}\nolimits}
\newcommand{\Rat}{\mathop{\text{\rm Rat}}\nolimits}
\newcommand{\A}{\mathop{\text{\rm A}}\nolimits}
\newcommand{\N}{\mathop{\text{\rm N}}\nolimits}

\newcommand{\Gal}{\mathop{\text{\rm Gal}}}
\newcommand{\Hom}{\mathop{\text{\rm Hom}}}
\newcommand{\summ}{\mathop{\text{\rm sum}}}
\newcommand{\ev}{\mathop{\text{\rm ev}}\nolimits}

\newcommand{\reg}{{\text{\rm reg}}}
\renewcommand{\top}{{\text{\rm top}}}

\newcommand{\calV}{\mathscr{V}}

\newcommand{\bbF}{{\mathbbm F}}
\newcommand{\bbG}{{\mathbbm G}}
\newcommand{\bbQ}{{\mathbbm Q}}
\newcommand{\bbZ}{{\mathbbm Z}}

\newcommand{\br}{ }
\newcommand{\brr}{, }

\def\rightend#1#2{{%
 \leavevmode\nobreak\hskip .5em plus 1fil
 \penalty600 \hskip 0pt plus -1filll
 \vadjust{}\nobreak\hskip 0pt plus 1filll%
 #1\parfillskip=#2\relax \par}}

\def\eop{\ifmmode\rule[-22pt]{0pt}{1pt}\ifinner\tag*{$\square$}\else\eqno{\square}\fi\else\rightend{$\square$}{0pt}\fi}

\renewcommand{\thefootnote}{\fnsymbol{footnote}}
\author{Andreas-Stephan Elsenhans${}^*$ and J\"org Jahnel${}^\ddagger$}
\date{}

\title{On the quasi group of a cubic surface\\ over a finite field}

\begin{document}
\renewcommand{\thefootnote}{\fnsymbol{footnote}}

\maketitle

\begin{abstract}
We construct nontrivial homomorphisms from the quasi group of some cubic surfaces
over~$\bbF_{\!p}$
into a~group. We~show experimentally that the homomorphisms constructed are the only possible ones and that there are no nontrivial homomorphisms in the other~cases. Thereby,~we follow the classification of cubic surfaces, due to A.\,Cayley. 
\end{abstract}

\footnotetext[1]{Mathematisches Institut, Universit\"at Bayreuth, Univ'stra\ss e 30, D-95440 Bayreuth, Germany,\\
{\tt Stephan.Elsenhans@uni-bayreuth.de}, Website: {\tt http://www.staff.uni-bayreuth.de/$\sim$btm216}}
\footnotetext[3]{Fachbereich 6 Mathematik, Universit\"at Siegen, Walter-Flex-Str.~3, D-57068 Siegen, Germany,\\
{\tt jahnel@mathematik.uni-siegen.de}, Website: {\tt http://www.uni-math.gwdg.de/jahnel}}
\footnotetext[1]{The first author was supported in part by the Deutsche Forschungsgemeinschaft (DFG) through a funded research~project.}
\footnotetext[1]{${}^\ddagger$Key words: Cubic surface, Cayley's classification, quasi group, Mordell-Weil group, Suslin's homology. MSC: 14J26, 14J20, 14G15, 11G25}

\section{The quasi group of a cubic surface}

\begin{ttts}
According~to Yu.\,I.~Manin~\cite{Ma}, a cubic
surface~$V$
carries a structure of a {\em quasi~group}. For~us, this shall simply mean the ternary relation
$$[x_1, x_2, x_3] \Longleftrightarrow x_1, x_2, x_3 \text{~non-singular,~intersection~of~}
V {\rm ~with~a~line} \, .$$
If~$V$
is defined over a
field~$K$
then,
on~$V^\reg (K)$,
there is a structure of a quasi~group.
\end{ttts}

\begin{ttts}
Here,~the precise definition is that the lines lying entirely on the surface shall {\em not\/} cause any~relation. On~the other hand, it is allowed that two or all three points~coincide. Then,~the line shall simply be tangent to the surface of order two or~three.
\end{ttts}

\begin{defis}
Let~$(\Gamma, [\,])$
be a quasi~group and
$(G, +)$
be an abelian~group. By~a homomorphism
$p\colon \Gamma \to G$,
we mean a mapping such that, for a
suitable~$g \in G$,
$$p(x_1) + p(x_2) + p(x_3) = g$$
whenever~$[x_1, x_2, x_3]$
is~true.
\end{defis}

\begin{facs}
The~category of all homomorphisms from a
quasi~group\/~$\Gamma$
to abelian groups has an initial~object. The~corresponding abelian group is\/
$\underline\Gamma := \bbZ\Gamma / N$,
for\/
$N$
the subgroup generated by\/
$1 \!\cdot\! x_1 + 1 \!\cdot\! x_2 + 1 \!\cdot\! x_3 - 1 \!\cdot\! x_1^\prime - 1 \!\cdot\! x_2^\prime - 1 \!\cdot\! x_3^\prime$
for all\/
$[x_1, x_2, x_3]$
and\/~$[x_1^\prime, x_2^\prime, x_3^\prime]$.\smallskip

\noindent
$\underline\Gamma$~carries
a surjective augmentation
homomorphism\/~$s \colon \underline\Gamma \to \bbZ$.
We~will call\/
$\ker s$
the\/ {\em group associated}
with\/~$\Gamma$.
\end{facs}

\begin{defiss}
Let~$V$
be a cubic surface over a
field~$K$.

\begin{iii}
\item
We~will call the group associated with the quasi~group
$V^\reg (K)$
the {\em Mordell-Weil group\/}
of~$V$.
It~will be denoted
by~$\MW(V)$.
\item
We~will call two points
$x_1, x_2 \in V^\reg(K)$
{\em equivalent\/} if
$[x_1] - [x_2] = 0 \in \MW(V)$.
\end{iii}
\end{defiss}

\begin{exes}
Let~$V$
be the {\em Cayley cubic\/} given by the equation
$$xyz + xyw + xzw + yzw = 0$$
in~$\Pb^3$
over a
field~$K$.
Then,~for a non-singular
point
$p = (x_0:y_0:z_0:w_0) \in V(K)$,
either no coordinate vanishes or exactly two of~them. Accordingly,~put
$$
s(p) := \left\{
\begin{array}{cl}
x_0y_0z_0w_0 & {\rm ~if~} x_0,y_0,z_0,w_0 \neq 0 \, ,\\
-\!\!\!\!\prod\limits_{v = x_0,y_0,z_0,w_0 \neq 0}\!\!\!\!\!\!\!\!\!\!\!\! v & {\rm ~otherwise} \, .
\end{array}
\right.
$$
Further,~let
$l$
be a line
in~$\Pb^3$
not contained
in~$V$
and denote by
$p_1$,
$p_2$,
and~$p_3$
the intersection points
with~$V$,
which are supposed to be non-singular and
\mbox{$K$-rational}
and counted with~multiplicity.\smallskip

\noindent
Then,~$s(p_1)s(p_2)s(p_3)$
is a perfect square
in~$K$.\medskip

\noindent
{\bf Proof.}
This~observation is easily checked by calculations in~{\tt maple}, treating the possible cases separately.
\eop
\end{exes}

\begin{exess}[{\rm continued}{}]
For~$V$
the Cayley cubic
over~$K$,
the map
$s$
therefore induces a surjective homomorphism of~groups
$$\MW(V) \longrightarrow K^*/(K^*) {}^2 \, .$$

\begin{iii}
\item
In~particular, for the Cayley cubic
over~$\bbQ$,
the group
$\MW(V)$
is not finitely~generated. 
\item
On~the other hand, for the Cayley cubic over a finite field
$\bbF_{\!q}$
of odd characteristic, we have
$\MW(V) \cong \bbZ/2\bbZ$.
There~are two different kinds of smooth points
on~$V$.
Two~points
$p_1, p_2 \in V^\reg(\bbF_{\!q})$
are equivalent if and only if
$s(p_1)s(p_2)$
is a~square.

The~purpose of this article is to investigate this phenomenon more~systematically.
\end{iii}
\end{exess}

\begin{remes}
The~Mordell-Weil group is related to the famous {\em Mordell-Weil problem,} which may be formulated as to find a minimal system of generators
for~$\MW(V)$.
\end{remes}

\begin{ttts}
We~are particularly interested in the cases when,
for~$V$
a cubic surface over a finite field,
$\MW(V) \neq 0$.
The~point is that there is the~following\smallskip

\noindent
{\bf Application.}
Let~$\calV$
be a cubic surface
over~$\bbQ$
and
$p_1, \ldots, p_t$
be primes satisfying the following~conditions.

\begin{iii}
\item
There~is no
\mbox{$\bbF_{\!p_i}$-rational}
line contained in the
reduction~$\calV_{p_i}$.
\item
The singularities of
$\calV_{p_i}$
do not lift to smooth
\mbox{$\bbQ$-rational}
points
on~$\calV$.
\end{iii}
Then,~the reductions induce a natural~homomorphism
$$\MW(\calV) \longrightarrow \MW(\calV_{p_1}) \times \ldots \times \MW(\calV_{p_t}) \, .$$
If~$\calV$
has weak approximation then this map is a~surjection.
\end{ttts}

\begin{ttts}
The~is another application, which is related to the so-called Brauer-Manin~obstruction. This~is a method, invented by Yu.\,I.~Manin~\cite[Chapter.~VI]{Ma}, to explain the failure of the Hasse principle or weak approximation in certain~cases. It~is based on the consideration of a non-trivial Brauer~class
$\alpha \in \Br(\calV)$
and the corresponding
\mbox{$p$-adic}
evaluation~maps
$$\ev_{\alpha,p} \colon \calV(\bbQ_p) \longrightarrow \Br(\bbQ_p) = \bbQ/\bbZ, \quad x \mapsto \alpha|_x \, .$$
{\bf Proposition.}
{\em
Let\/~$p \neq 2$
be a prime number and\/
$F \in \bbZ_p[X_0,X_1,X_2,X_3]$
cubic form defining a smooth cubic surface\/
$\calV$
over\/~$\bbQ_p$.
Suppose~that all\/
$x \in \calV(\bbQ_p)$
specialize
to\/~$\calV_p^\reg$
and that\/
$\MW(\calV_p) = 0$.\smallskip

\noindent
Then,~$\ev_{\alpha,p}$
is constant for
every\/~$\alpha \in \Br(\calV)$.}\medskip

\noindent
{\bf Proof.}
It is known that
$\ev_{\alpha,p}(x)$
depends only on the reduction of
$x$
modulo~$p$\break
\cite[Theorem~1]{B}.
Further,~an application of Lichtenbaum duality~\cite[Corollary~1]{Li} proves that
$\ev_{\alpha,p}$
is induced by a group homomorphism
$\MW(\calV_p) \to \bbQ/\bbZ$.
\eop
\end{ttts}

\begin{remes}
In~\cite{EJ}, we studied explicit examples of cubic surfaces, for which the Brauer-Manin obstruction works at certain~primes. It~was noticeable in the experiments that the reduction types at the relevant primes were distributed in an unusual~way. Reducible~reductions and reductions to the Cayley~cubic occurred~frequently. This~observation was actually the starting point of our investigations on the Mordell-Weil~group.
\end{remes}

\begin{ttts}
The~goal of this article is to compare
$\MW(V)$
with a group more tractable from the theoretical point of~view. For~cubic surfaces that are not too degenerate, this will be
$\A_0(V^\reg)$,
the
degree-$0$~part
of Suslin's homology
group~$h_0(V^\reg)$.
We~will establish a canonical homomorphism
$$\pi_V \colon \MW(V) \longrightarrow \A_0(V^\reg)$$
for
$V$
a geometrically irreducible cubic surface over a finite~field. Under~minimal assumptions,
$\pi_V$
will be~surjective.
\end{ttts}

\begin{pls}
In~section~\ref{Cayley}, we will recall Cayley's classification of cubic~surfaces. After~this, we will consider two degenerate cases at~first. Section~\ref{crvs} will be concerned with the situation of a~cone. Then,~there is a surjection to the Mordell-Weil group of the underlying~curve. Section~\ref{red} will treat the reducible~case. It~will turn out that there is a nontrivial surjection
from~$\MW(V)$
to a nontrivial abelian group, which is given in an elementary~manner. In~section~\ref{Suslin}, we will construct the
homomorphism~$\pi_V$.
Then,~we will compute 
$\A_0(V^\reg)$
systematically for each of the remaining cases of the classification of cubic~surfaces.

At~the end of the main body of the article, we will report on the comparison between
$\MW(V)$
and~$\A_0(V^\reg)$
in a large sample of~examples.
In~an appendix, we will discuss efficient algorithms to compute
$\MW(V)$
for a concrete~surface.
\end{pls}

\section{Cayley's classification of cubic surfaces}
\label{Cayley}

\begin{ttts}
Cubic~surfaces are classified since the days of A.\,Cayley~\cite[sec.\,9.2]{Do}. According~to this, there are the following~types.

\begin{III}
\item
A~normal cubic surface is either
\begin{iii}
\item
in one of the 21~classes of surfaces with finitely many double points, listed in~\cite[Table~9.2.5]{Do}. This~includes the case of a smooth cubic~surface.
\item
Or~the cone over a smooth cubic
curve~$C$.
\end{iii}
\item
A non-normal, geometrically irreducible cubic surface is either
\begin{iii}
\item
a cubic ruled~surface. There~are two types of those~\cite[Theorem~9.2.1]{Do}, ordinary and Cayley's cubic~ruled~surfaces.
\item
Or~the cone over a singular cubic~curve. This~might be a cubic with a self-inter\-section or a~cusp.
\end{iii}
\end{III}
\end{ttts}

\begin{ttts}
\label{classif}
In~the situation of a finite base field, the classification of geometrically irreducible cubic surfaces is actually a little~finer.

\begin{iii}
\item[I.i) ]
Among these types,
$2A_1$,
$3A_1$,
$2A_2$,
$A_2 + 2A_1$,
$4A_1$,
$2A_2 + A_1$,
$A_3 + 2A_1$,
and~$3A_2$
have~symmetries. This~leads to 13 further types, where the singularities are defined over extensions of the ground~field.
\item[II.i) ]
An~ordinary cubic ruled surface may have its normal form
$xz^2 + yw^2 = 0$
only over a quadratic~extension. This~causes a third type of cubic ruled surfaces over a finite~field.
\item[II.ii) ]
In~the case of the cone over a cubic curve with self-inter\-section, there are two variants as to whether the two tangent directions at the point of intersection are defined over the ground field or~not.
\end{iii}
\end{ttts}

\begin{ttts}
\label{classifred}
We~will restrict ourselves to reduced cubic~surfaces. In~other words, the following types of reducible surfaces are~allowed.

\begin{iii}
\item
A~reducible cubic surface might consist of a quadric and a~plane. There~are four cases where the quadric is nondegenerate. In~fact, the quadric may split over the ground field or not and the plane may be tangent or~not. There~are four more cases when the quadric is a~cone. The~intersection with the plane might be a conic, two lines, a double line, or a point.
\item
Finally,~the surface might be reducible into three~planes. There~are two cases as to whether their intersection is a point or a~line. Observe~that it is possible that the decomposition into three planes is defined only after a finite field~extension.
\end{iii}
\end{ttts}

\section{The case of a curve}
\label{crvs}

\begin{ttts}
Let~$C \subset \Pb^2$
be a reduced cubic curve over a
field~$K$.
Then,~in a manner analogous to the surface case, there are a quasi group structure
on~$C$
and the
{\em Mordell-Weil group\/}~$\MW(C)$.
This~group is known in every~case.

\begin{iii}
\item
It~may happen that all
\mbox{$K$-rational}
smooth points are contained in a~line. Then,~the quasi~group structure is empty
and~$\MW(C) = \ker(\summ\colon \bbZ[L(K)] \to \bbZ)$.
We~have~this degenerate case whenever
$C$
contains a line defined over a proper extension
of~$K$.
The~same may happen even for a smooth cubic curve when
$\#K \leq 5$.

Otherwise,
\item
$\MW(C) = J(C)(K)$
for
$C$
smooth, and
\item
$\MW(C) = K^+$
if
$C$
is a cubic curve with a~cusp.
\item
If~$C$
is a cubic curve with a node~then
$\MW(C) = K^*$
in case that the two tangent directions at the node are defined
over~$K$.
If~the tangent directions are defined over the quadratic
extension~$F/K$
then
$\MW(C) = \ker (\N\colon F^* \to K^*)$.
\item
When~$C$
is reducible into a line and a conic then
$\MW(C) = \ker (\N\colon F^* \to K^*) \oplus \bbZ$,
$\MW(C) = K^+ \oplus \bbZ$
or
$\MW(C) = K^* \oplus \bbZ$
depending on whether,
over~$K$,
there are no, one, or two points of intersection.
\item
When~$C$
is reducible into three components then
$\MW(C) = K^+ \oplus \bbZ^2$
or
$\MW(C) = K^* \oplus \bbZ^2$
depending on whether the three points of intersection coincide or~not. In~the case of three points of intersection, this is actually Menelaos'~Theorem.
\end{iii}
\end{ttts}

\begin{facos}[{\rm Cones}{}]
For\/~$V$
a cone over a cubic
curve\/~$C$,
we have a canonical
surjection\/~$\MW(V) \to \MW(C)$.
\end{facos}

\section{Reducible cubic surfaces}
\label{red}

\begin{ttts}
\label{redu}
Let~$V$
be a reducible cubic surface over a
field~$K$.
Then,~there are two essentially different~cases.

\begin{iii}
\item
There~are two irreducible components, a
plane~$E$
and a quadric, but the quadric consists of two planes defined over a quadratic~extension. Then,~only the
plane~$E$
contains
\mbox{$K$-rational}
smooth~points. We~have an empty quasi group structure
and~$\MW(V) = \ker(\summ\colon \bbZ[E(K)] \to \bbZ)$.
\item
Otherwise,~when
$V$
decomposes
into~$k = 2,3$
components, there is a canonical~surjection
$\MW(V) \twoheadrightarrow \ker(\summ\colon \bbZ^k \to \bbZ) \cong \bbZ^{k-1}$.
\end{iii}
\end{ttts}

\begin{exes}
\label{4klassen}
Over~a finite
field~$\bbF_{\!q}$
of characteristic
$\neq\!2$,
let~$V$
be a reducible cubic surface consisting of a nondegenerate quadratic
cone~$Q$
and a
plane~$E$.
Suppose~that
$E$
does not meet the cusp
of~$Q$.
Then,~there is a canonical~surjection
$$\MW(V) \twoheadrightarrow \bbZ \oplus \bbZ/2\bbZ \, .$$
{\bf Proof.}
The~homomorphism
to~$\bbZ$
is that from~\ref{redu}.ii). It~remains to construct the homomorphism
to~$\bbZ/2\bbZ$.

For~this, we fix coordinates such that the cusp is
in~$(1:0:0:0)$
and the
plane~$E$
is given
by~$x=0$.
Further,~we assume without restriction that the plane
``$y=0$''
is tangent to the
cone~$Q$.
Then,~the cone is given by,
say,~$yz+Kw^2 = 0$
for~$K \neq 0$.
The~whole cubic surface has the~equation
$$x(yz+Kw^2) = 0 \, .$$
On~the
plane~``$x=0$'',
we define the homomorphism
$\MW(V) \to \bbZ/2\bbZ$
simply as
$\chi_2\big(\!K(yz+Kw^2)\big)$
for
$\chi_2$
the quadratic character
on~$\bbF_{\!q}^*$.
On~the
cone~``$yz+Kw^2 = 0$''\!,
we
take~$\chi_2(xy)$,
respectively
$\chi_2(-Kxz)$
when~$y=0$.

We~have to show that this definition is indeed compatible with the quasi group~structure. For~this, let
$(x:y:z:w) \in Q(\bbF_{\!q})$,
$(x':y':z':w') \in Q(\bbF_{\!q})$,
and~$(0:y'':z'':w'') \in E(\bbF_{\!q})$
be three collinear~points. Then,~we clearly have
$(0:y'':z'':w'') = (0 : (x'y - xy') : (x'z - xz') : (x'w - xw'))$.
Furthermore,
\begin{eqnarray*}
(xy)(x'y')K(y''z'' + Kw''^2) & = & Kxx'yy'[(x'y - xy')(x'z - xz') + K(x'w - xw')^2] \\
 & = & -Kxx'yy'[xx' (yz' + y'z + 2Kww')] \\
 & = & -Kx^2 x'^2 (-Ky^2w'^2 - Ky'^2w^2 + 2Kyy'ww') \\
 & = & K^2x^2 x'^2 (yw' - y'w)^2
\end{eqnarray*}
is a perfect~square.
\eop
\end{exes}

\section{Irreducible cubic surfaces not being cones}
\label{Suslin}

\begin{ttt}
When~a cubic
surface~$V$
is irreducible, but geometrically reducible, then it consists of three planes acted upon transitively by the Galois~group. In~this case,
$V^\reg(K) = \emptyset$
and,
therefore,~$\MW(V) = 0$.
Thus,~we may restrict ourselves to the geometrically irreducible~case.
\end{ttt}

\subsection{\boldmath Suslin's singular homology
group~$h_0$}

\begin{ttt}
For~a scheme of finite type over a
field~$K$,
the singular homology groups
$h_*(S)$
were introduced by A.~Suslin~\cite{SV}. We~will only need 
$h_0(S)$,
for which there is the following elementary~description.
\end{ttt}

\begin{fac}
Let\/~$S$
be an integral scheme of finite type over a
field\/~$K$.
Then,
$$h_0(S) = \Z_0(S) / \Rat_0^\prime(S) \, .$$
Here,
$\Z_0 (S)$
is the group of\/
\mbox{$0$-cycles},
i.e., the free abelian group over all closed points
of\/~$S$.
$\Rat_0^\prime(S)$~is
generated by all\/
\mbox{$0$-cycles}
of the following~kind.\smallskip

\noindent
Let\/~$C \subset S$
be an irreducible curve,
$C^\prime$
its normalization, and\/
$\overline{C}$
the corresponding smooth, proper~model. Then,~take all the cycles\/
$\div(f)$
where\/
$f \not\equiv 0$
is a rational function\/
on\/~$C$
which, after pull-back
to\/~$\overline{C}$,
is
constantly\/~$1$
on\/~$\overline{C} \!\setminus\! C^\prime$.\smallskip

\noindent
{\bf Proof.}
{\em
See~\cite[Theorem~5.1]{Sch}.
}
\eop
\end{fac}

\begin{rems}
\begin{abc}
\item
$h_0 (S)$~is
equipped with a natural map
$\deg\colon h_0 (S) \to \bbZ$.
We~will denote its kernel
by~$\A_0 (S)$.
\item
Let~$i\colon S_1 \to S_2$
be an arbitrary morphism of quasi-projective varieties
over~$K$.
Then,~there is the induced~homomorphism
$i_* \colon h_0 (S_1) \longrightarrow h_0 (S_2), [x] \mapsto [i(x)]$.
This~immediately yields a
map~$i_* \colon \A_0 (S_1) \longrightarrow \A_0 (S_2)$.
\end{abc}
\end{rems}

\begin{lem}
Let\/~$V$
be a geometrically irreducible cubic surface
over\/~$\bbF_{\!q}$.
Then,~there is a canonical~homomorphism
$$\pi_V \colon \MW(V) \longrightarrow \A_0 (V^\reg) \, .$$
{\bf Proof.}
{\em
To~each combination
$a_1 [p_1] + \ldots + a_k [p_k]$
for
$p_1, \ldots, p_k \in V^\reg(K)$
and
$a_1 + \ldots + a_k = 0$,
the homomorphism
$i_*$
assigns the corresponding~cycle. We~take this as a definition
for~$\pi_V$.
To~show that
$\pi_V$
is well-defined, we have to verify the~following.

Assume~that
$x_1, x_2, x_3$
are collinear and
$x_1^\prime, x_2^\prime, x_3^\prime$
are collinear,~too. Suppose~that the connecting lines are not contained
in~$V$.
Then,
$$[x_1] + [x_2] + [x_3] - [x_1^\prime] - [x_2^\prime] - [x_3^\prime] = 0 \in \A_0(V^\reg) \, .$$

For~this, consider the pencil of planes
through~$x_1, x_2, x_3$.
Generically,~the intersection
with~$V$
is a curve, smooth at
$x_1$,
$x_2$
and~$x_3$.
The~only possible exceptions are the tangent~planes. We~claim that the generic intersection curve is irreducible,~too. Indeed,~the contrary would mean that all intersection curves contained a~line. Suppose,~this is a line
through~$x_1$.
Then,~$V$
contains a pencil of lines
through~$x_1$,
which implies
$V$
contains a plane
through~$x_1$.
Hence,~$V$
is reducible, a~contradiction.

Thus,~take a plane
through~$x_1, x_2, x_3$,
generating an irreducible intersection
curve~$C$
that is smooth in
$x_1$,
$x_2$
and~$x_3$.
Further,~take a plane
through~$x_1^\prime, x_2^\prime, x_3^\prime$
generating an irreducible intersection
curve~$C^\prime$
that is smooth in
$x_1^\prime$,
$x_2^\prime$
and~$x_3^\prime$
and meets
$C$
only in smooth
points~$x_1^{\prime\prime}, x_2^{\prime\prime}, x_3^{\prime\prime}$.
The~sublemma below, applied
to~$C$
and~$C^\prime$,
immediately yields the~assertion.
}
\eop
\end{lem}

\begin{sub}
Let\/~$C$
be an irreducible cubic~curve. Assume~that\/
$p_1, p_2, p_3 \in C^\reg$
as well as\/
$q_1, q_2, q_3 \in C^\reg$
are triples of collinear points such that
$\{p_1, p_2, p_3\} \cap \{q_1, q_2, q_3\} = \emptyset$.\smallskip

\noindent
Then,~there is a rational function\/
$f$
on\/~$C$
having simple zeroes
at~\/$p_1, p_2, p_3$,
simple poles
at~\/$q_1, q_2, q_3$,
no other zeroes or poles, and the value\/
$1$
at the possible singular~point.\smallskip

\noindent
{\bf Proof.}
{\em
According~to J.\,Pl\"ucker, an irreducible cubic curve may have at most one singular~point. We~may therefore
put~$f := K \!\cdot\! l_1/l_2$
for forms
$l_1$
and~$l_2$
defining the~lines. By~assumption, these do not meet the singular~point. If~necessary, we choose the
constant~$K$
such that the value at the singularity is normalized
to~$1$.
}
\eop
\end{sub}

\subsection{\boldmath
$h_0$~and
the tame fundamental group}

\begin{ttt}
Let~$S$
be a smooth surface over the finite
field~$\bbF_{\!q}$
for~$q = p^r$
and let
$\overline{S} \supseteq S$
be a smooth~compactification. Then,~the {\em tame fundamental group\/}
$\pi_1^t(S)$
of~$S$
classifies all finite coverings
of~$S$
which are tamely ramified
at~$\overline{S} \!\setminus\! S$.

The~group
$\pi_1^t(S)$
is independent of the choice of the
compactification~$\overline{S}$.
$\pi_1^t(S)$~is
a quotient
of~$\pi_1^\et(S)$.
By~the purity of the branch locus \cite[Exp.~X, Th\'eor\`eme~3.1]{SGA1}, one~has
$$\pi_1^t(S)_\tors^\ab \cong (\pi_1^\et(S)^\ab)_{\text{prime\,to\,}p} \oplus (\pi_1^\et(\overline{S})^\ab)_{p\text{-power}}\,.$$
Again,~this decomposition is independent of the choice
of~$\overline{S}$.

The~structural morphism
$S \to \Spec \bbF_{\!q}$
induces a
surjection~$\pi_1^t(S) \to \pi_1(\Spec \bbF_{\!q})$
the kernel of which we will denote
by~$\pi_1^{t,\geo}(S)$.
Note~that
$\pi_1^{t,\geo}(S)$
differs from
$\pi_1^t(S_{\overline\bbF_{\!q}})$.
The~point is that the analogue of the natural short exact sequence \cite[Exp.~IX, Th\'eor\`eme~6.1]{SGA1} is only right exact for the tame fundamental~group.
\end{ttt}

\begin{theoremo}[{\rm Schmidt, Spie{\ss}}{}]
Let\/~$S$
be a surface over a finite
field\/~$\bbF_{\!q}$
which is smooth and geometrically irreducible, but not necessarily~proper.

\begin{iii}
\item
Then,~$\A_0 (S)$
is a finite abelian~group.
\item
There~is a canonical~isomorphism\/
$\iota_S \colon \A_0 (S) \longrightarrow \pi_1^{t,\geo}(S)^\ab$.
\end{iii}
{\bf Proof.}
{\em
See~\cite[Theorem~0.1]{SchS}.
\eop
}
\end{theoremo}

\begin{rems}
\begin{abc}
\item
Concretely,~$\iota_S$
is given as~follows.
\begin{iii}
\item
For~a point
$x\colon \Spec \bbF_{\!q'} \rightarrow S$,
consider the induced homomorphism
$$\pi_1^\et(x) \colon \widehat{\bbZ} = \pi_1^\et(\Spec \bbF_{\!q'}\!) \longrightarrow \pi_1^\et(S) \twoheadrightarrow \pi_1^t(S) \twoheadrightarrow \pi_1^t(S)^\ab \, .$$
Send~$[x]$
to~$\pi_1^\et(x)(1)$.
This~defines a homomorphism
$\iota_S^\prime \colon h_0 (S) \to \pi_1^t(S)^\ab$.
\item
Clearly,~the degree map
$\deg \colon h_0 (S) \to \bbZ$
is compatible with the homomorphism
$\smash{\pi_1^t(S)^\ab \to \pi_1^\et(\Spec \bbF_{\!q}) = \widehat{\bbZ}}$
induced by the structural~morphism.
\item
The~homomorphism
$\iota_S$
is exactly the restriction
of~$\smash{\iota_S^\prime}$
to~$\ker\deg$.
\end{iii}
\item
The~map
$\iota_S^\prime$
defines an isomorphism
$\smash{\widehat{h_0 (S)} \to \pi_1^t(S)^\ab}$.
\end{abc}
\end{rems}

\subsection{The tame fundamental group and the Picard~group}

\begin{fac}
Let\/~$V$
be a cubic ruled surface defined over the finite
field\/~$\bbF_{\!q}$.
Then,
$\pi_1^{t,\geo} (V^\reg) = 0$.\smallskip

\noindent
{\bf Proof.}
{\em
It~will suffice to
show~$\smash{\pi_1^t (V^\reg_{\overline\bbF_{\!q}}) = 0}$.
In~the present situation, a smooth compactification
of~$\smash{V^\reg_{\overline\bbF_{\!q}}}$
is given by a projective plane, blown up in one~point. The~preimage of the singular locus is a (double) line through the point blown~up.
Consequently,~$\smash{V^\reg_{\overline\bbF_{\!q}}}$
is a ruled surface
over~$\Ab^1$.
This~yields~$\smash{\pi_1^t (V^\reg_{\overline\bbF_{\!q}}) = 0}$.
}
\eop
\end{fac}

\begin{prop}
Let\/~$V$
be a geometrically irreducible cubic surface
over\/~$\bbF_{\!q}$
that is not a~cone.
Suppose\/~$V$
is normal, i.e., of one of the types~I.i). Then,
$$\pi_1^{t,\geo}(V^\reg)^\ab = [(\Pic(V^\reg)_{{\rm prime\,to\,}p} \otimes_\bbZ \mu_\infty^\vee)^{\Gal(\overline\bbF_{\!q}/\bbF_{\!q})}]^\vee \, .$$
Here,~${}^\vee$
denotes the Pontryagin~dual, given by the
functor\/~$\Hom(\cdot, \bbQ/\bbZ)$.\smallskip\pagebreak[3]

\noindent
{\bf Proof.}
{\em
{\em First step.}
$p$-torsion.\smallskip

\noindent
We~know a smooth
compactification~$\overline{V}$
of~$V^\reg$,~explicitly.
$\smash{\overline{V}_{\!\overline\bbF_{\!q}}}$~is
isomorphic to
$\Pb^2$
blown-up in six~points. In~particular, we
have~$\smash{\pi_1^\et(\overline{V}_{\!\overline\bbF_{\!q}}) = 0}$.
This~suffices for
$\smash{\pi_1^t(V_{\overline\bbF_{\!q}})^\ab_{p\text{-power}} = 0}$
and~$\smash{\pi_1^{t,\geo}(V)^\ab_{p\text{-power}} = 0}$.\medskip

\noindent
{\em Second step.}
The Pontryagin dual.\smallskip

\noindent
Let~us compute the Pontryagin
dual~$(\pi_1^{t,\geo}(V^\reg)^\ab)^\vee$.
For~$l$
prime
to~$p$,
we~have
\begin{eqnarray*}
\textstyle \Hom(\pi_1^{t,\geo}(V^\reg)^\ab, \frac1l\bbZ/\bbZ)
 & = & \textstyle \Hom(\pi_1^t(V^\reg)^\ab, \frac1l\bbZ/\bbZ) / \Hom(\pi_1(\Spec \bbF_{\!q}), \frac1l\bbZ/\bbZ) \\
 & = & \textstyle \Hom(\pi_1(V^\reg)^\ab, \frac1l\bbZ/\bbZ) / \Hom(\pi_1(\Spec \bbF_{\!q}), \frac1l\bbZ/\bbZ) \\
 & = & \smash{\textstyle H^1_\et (V^\reg, \frac1l\bbZ/\bbZ) / H^1(\Gal(\overline\bbF_{\!q}/\bbF_{\!q}), \frac1l\bbZ/\bbZ) \, .}
\end{eqnarray*}
According~to the Hoch\-schild-Serre spectral~sequence
$$\smash{\textstyle H^p(\Gal(\overline\bbF_{\!q}/\bbF_{\!q}), H^q_\et (V^\reg_{\overline\bbF_{\!q}}, \frac1l\bbZ/\bbZ)) \Longrightarrow H^{p+q}_\et (V^\reg, \frac1l\bbZ/\bbZ) \, ,}$$
the latter quotient is nothing but
$\smash{H^1_\et (V^\reg_{\overline\bbF_{\!q}}, \frac1l\bbZ/\bbZ)^{\Gal(\overline\bbF_{\!q}/\bbF_{\!q})}}$.\medskip

\noindent
{\em Third step.}
The torsion part of the Picard group.\smallskip

\noindent
We~have
$\smash{\Gamma(V^\reg_{\overline\bbF_{\!q}}, \bbG_m) = \overline\bbF_{\!q}^*}$.
In~fact, the
$A$-,
$D$-,
and
$E$-configurations
do not contain any principal~divisor.
This~immediately yields
$\smash{H^1_\et (V^\reg_{\overline\bbF_{\!q}}, \mu_l) = \Pic(V^\reg_{\overline\bbF_{\!q}})_l}$
for any
$l$
prime
to~$p$.
On~$\smash{V^\reg_{\overline\bbF_{\!q}}}$,
the sheaves
$\mu_l$
and~$\frac1l\bbZ/\bbZ$
coincide up to the Galois~operation. We~therefore have
\begin{eqnarray*}
\textstyle \Hom(\pi_1^{t,\geo}(V^\reg)^\ab, \frac1l\bbZ/\bbZ) & = & (H^1_\et (V^\reg_{\overline\bbF_{\!q}}, \mu_l) \otimes_\bbZ \mu_l^\vee)^{\Gal(\overline\bbF_{\!q}/\bbF_{\!q})} \\
 & = & (\Pic(V^\reg_{\overline\bbF_{\!q}})_l \otimes_\bbZ \mu_l^\vee)^{\Gal(\overline\bbF_{\!q}/\bbF_{\!q})} \, .
\end{eqnarray*}
Summing this up over
all~$l$,
we see that
$$(\pi_1^{t,\geo}(V^\reg)^\ab)^\vee = (\Pic(V^\reg)_{{\rm prime\,to\,}p} \otimes_\bbZ \mu_\infty^\vee)^{\Gal(\overline\bbF_{\!q}/\bbF_{\!q})} \, ,$$
which is equivalent to the~assertion.
}
\eop
\end{prop}

\begin{summa}
Thus,~in order to compute
$\A_0(V^\reg)$,
we only need to know
$\Pic(V^\reg)_\tors$
for each of the 21 types of cubic surfaces summarized in~I.i).
\end{summa}

\subsection{The 21 types of normal cubic surfaces not being cones}

\begin{fac}
\label{quot}
Let\/
$V$
be a normal, proper surface over an algebraically closed field and\/
$\widetilde{V}$
its~desingularization. Then,
$$\Pic(V^\reg) = \Pic(\widetilde{V}) / \langle E_1, \ldots, E_k \rangle \, ,$$
where\/
$E_1, \ldots, E_k$
denote the irreducible components of the preimages of the singularities
on\/~$\widetilde{V}$.
\eop
\end{fac}

\begin{theorem}
Let\/
$V$
be an irreducible cubic surface
over\/~$\overline{\bbF}_{\!q}$,
not being a~cone.
Suppose~that\/
$V$
is normal, i.e., of one of the 21~types~I.i).\smallskip

\noindent
Then,~the Picard~group\/
$\Pic(V^\reg)$
is torsion-free for 17 of the 21~types. For~the four remaining types, the torsion is given in the table~below.

\begin{table}[H]
\scriptsize
\begin{center}
\begin{tabular}{|c|c|c|}
\hline
type & singularities & $\Pic(V^\reg)_\tors$ \\\hline
\hline
XVI & $4A_1$ & $\bbZ/2\bbZ$ \\\hline
XVIII & $A_3 + 2A_1$ & $\bbZ/2\bbZ$ \\\hline
XIX & $A_5 + A_1$ & $\bbZ/2\bbZ$ \\\hline
XXI & $3A_2$ & $\bbZ/3\bbZ$ \\\hline
\end{tabular}
\end{center}
\normalsize
\end{table}\vskip-0.4cm

\noindent
{\bf Proof.}
{\em
We~distinguish the cases~systematically. Each~time, we apply Fact~\ref{quot}.\smallskip

\noindent
One~has
$\Pic(\widetilde{V}) \cong \bbZ^7$.
The~signature
is~$(1,-1,-1,-1,-1,-1,-1)$.
I.e.,~we have torsion-freeness in the case of a smooth cubic~surface.\smallskip

\noindent
Otherwise,~the
\mbox{$A$-,}
\mbox{$D$-,}
or
$E$-configuration
of~$(-2)$-curves
generates a sublattice
of~$\smash{\Pic(\widetilde{V})}$.
The~quotient has torsion if and only if this sublattice can be refined
in~$\bbZ^7$
without enlarging the~rank. This~immediately shows torsion-freeness in the cases
$A_n$
for~$n \neq 3$
and~$E_6$
as the lattice discriminants are~square-free.\smallskip

\noindent
For~the other cases, the constructions described in~\cite[page~278]{Do} yield explicit generators for sublattices
of~$\bbZ^7$.
We~summarize them in the following~table.

\begin{table}[H]
\scriptsize
\begin{center}
\begin{tabular}{|c|l|}
\hline
$2A_1$ &
\phantom{$1.~A_1\colon\,$}
$
\begin{array}{c}
(2,-1,-1,-1,-1,-1,-1) \\
(0,\phantom{-}0,\phantom{-}0,\phantom{-}0,\phantom{-}0,\phantom{-}1,-1)
\end{array}
$
\\\hline
$A_3$ &
\phantom{$1.~A_1\colon\,$}
$
\begin{array}{c}
(0,\phantom{-}0,\phantom{-}0,\phantom{-}0,\phantom{-}0,\phantom{-}1,-1) \\
(0,\phantom{-}0,\phantom{-}0,\phantom{-}0,\phantom{-}1,-1,\phantom{-}0) \\
(0,\phantom{-}0,\phantom{-}0,\phantom{-}1,-1,\phantom{-}0,\phantom{-}0)
\end{array}
$
\\\hline
$A_2 + A_1$ &
\phantom{$1.$}
$\begin{array}{ccrr}
A_1\colon & (2,-1,-1,-1,-1,-1,-1) && \phantom{=: v_1} \\
A_2\colon & (0,\phantom{-}0,\phantom{-}0,\phantom{-}0,\phantom{-}0,\phantom{-}1,-1) && \\
A_2\colon & (0,\phantom{-}0,\phantom{-}0,\phantom{-}0,\phantom{-}1,-1,\phantom{-}0) &&
\end{array}$
\\\hline
$3A_1$ &
\phantom{$1.~A_1\colon\,$}
$
\begin{array}{c}
(2,-1,-1,-1,-1,-1,-1) \\
(0,\phantom{-}0,\phantom{-}0,\phantom{-}0,\phantom{-}0,\phantom{-}1,-1) \\
(0,\phantom{-}0,\phantom{-}0,\phantom{-}1,-1,\phantom{-}0,\phantom{-}0)
\end{array}
$
\\\hline
$2A_2$ &
$
\begin{array}{rc}
1.~A_2\colon & (0,\phantom{-}0,\phantom{-}0,\phantom{-}0,\phantom{-}0,\phantom{-}1,-1) \\
1.~A_2\colon &
(0,\phantom{-}0,\phantom{-}0,\phantom{-}0,\phantom{-}1,-1,\phantom{-}0) \\
2.~A_2\colon &
(0,\phantom{-}0,\phantom{-}1,-1,\phantom{-}0,\phantom{-}0,\phantom{-}0) \\
2.~A_2\colon &
(0,\phantom{-}1,-1,\phantom{-}0,\phantom{-}0,\phantom{-}0,\phantom{-}0)
\end{array}
$
\\\hline
$A_3 + A_1$ &
\phantom{$1.$}
$\begin{array}{ccrr}
A_1\colon & (2,-1,-1,-1,-1,-1,-1) && \phantom{=: v_1} \\
A_3\colon & (0,\phantom{-}0,\phantom{-}0,\phantom{-}0,\phantom{-}0,\phantom{-}1,-1) && \\
A_3\colon & (0,\phantom{-}0,\phantom{-}0,\phantom{-}0,\phantom{-}1,-1,\phantom{-}0) && \\
A_3\colon & (0,\phantom{-}0,\phantom{-}0,\phantom{-}1,-1,\phantom{-}0,\phantom{-}0) &&
\end{array}$
\\\hline
$D_4$ &
\phantom{$1.~A_1\colon\,$}
$
\begin{array}{c}
(1,-1,-1,-1,\phantom{-}0,\phantom{-}0,\phantom{-}0) \\
(0,\phantom{-}1,\phantom{-}0,\phantom{-}0,-1,\phantom{-}0,\phantom{-}0) \\
(0,\phantom{-}0,\phantom{-}1,\phantom{-}0,\phantom{-}0,-1,\phantom{-}0) \\
(0,\phantom{-}0,\phantom{-}0,\phantom{-}1,\phantom{-}0,\phantom{-}0,-1)
\end{array}
$
\\\hline
$A_2 + 2A_1$ &
$\begin{array}{rc}
1.~A_1\colon & (2,-1,-1,-1,-1,-1,-1) \\
2.~A_1\colon & (0,\phantom{-}0,\phantom{-}0,\phantom{-}0,\phantom{-}0,\phantom{-}1,-1) \\
A_2\colon & (0,\phantom{-}0,\phantom{-}0,\phantom{-}1,-1,\phantom{-}0,\phantom{-}0) \\
A_2\colon & (0,\phantom{-}0,\phantom{-}1,-1,\phantom{-}0,\phantom{-}0,\phantom{-}0)
\end{array}$
\\\hline
$A_4 + A_1$ &
\phantom{$1.$}
$\begin{array}{ccrr}
A_1\colon & (2,-1,-1,-1,-1,-1,-1) && \phantom{=: v_1} \\
A_4\colon & (0,\phantom{-}0,\phantom{-}0,\phantom{-}0,\phantom{-}0,\phantom{-}1,-1) && \\
A_4\colon & (0,\phantom{-}0,\phantom{-}0,\phantom{-}0,\phantom{-}1,-1,\phantom{-}0) && \\
A_4\colon & (0,\phantom{-}0,\phantom{-}0,\phantom{-}1,-1,\phantom{-}0,\phantom{-}0) && \\
A_4\colon & (0,\phantom{-}0,\phantom{-}1,-1,\phantom{-}0,\phantom{-}0,\phantom{-}0) &&
\end{array}$
\\\hline
$D_5$ &
\phantom{$1.~A_1\colon\,$}
$
\begin{array}{c}
(1,-1,-1,\phantom{-}0,\phantom{-}0,\phantom{-}0,-1) \\
(0,\phantom{-}1,-1,\phantom{-}0,\phantom{-}0,\phantom{-}0,\phantom{-}0) \\
(0,\phantom{-}0,\phantom{-}1,-1,\phantom{-}0,\phantom{-}0,\phantom{-}0) \\
(0,\phantom{-}0,\phantom{-}0,\phantom{-}1,-1,\phantom{-}0,\phantom{-}0) \\
(0,\phantom{-}0,\phantom{-}0,\phantom{-}0,\phantom{-}1,-1,\phantom{-}0) 
\end{array}
$
\\\hline
$4A_1$ &
\phantom{$1.~A_1\colon\,$}
$
\begin{array}{c}
(2,-1,-1,-1,-1,-1,-1) \\
(0,\phantom{-}0,\phantom{-}0,\phantom{-}0,\phantom{-}0,\phantom{-}1,-1) \\
(0,\phantom{-}0,\phantom{-}0,\phantom{-}1,-1,\phantom{-}0,\phantom{-}0) \\
(0,\phantom{-}1,-1,\phantom{-}0,\phantom{-}0,\phantom{-}0,\phantom{-}0)
\end{array}
$
\\\hline
$2A_2 + A_1$ &
$
\begin{array}{rc}
A_1\colon & (2,-1,-1,-1,-1,-1,-1) \\
1.~A_2\colon & (0,\phantom{-}0,\phantom{-}0,\phantom{-}0,\phantom{-}0,\phantom{-}1,-1) \\
1.~A_2\colon & (0,\phantom{-}0,\phantom{-}0,\phantom{-}0,\phantom{-}1,-1,\phantom{-}0) \\
2.~A_2\colon & (0,\phantom{-}0,\phantom{-}1,-1,\phantom{-}0,\phantom{-}0,\phantom{-}0) \\
2.~A_2\colon & (0,\phantom{-}1,-1,\phantom{-}0,\phantom{-}0,\phantom{-}0,\phantom{-}0)
\end{array}
$
\\\hline
\end{tabular}
\end{center}
\normalsize
\end{table}

\begin{table}[H]
\scriptsize
\begin{center}
\begin{tabular}{|c|l|}
\hline
$A_3 + 2A_1$ &
$
\begin{array}{rc}
1.~A_1\colon & (2,-1,-1,-1,-1,-1,-1) \\
2.~A_1\colon & (0,\phantom{-}0,\phantom{-}0,\phantom{-}0,\phantom{-}0,\phantom{-}1,-1) \\
A_3\colon & (0,\phantom{-}0,\phantom{-}0,\phantom{-}1,-1,\phantom{-}0,\phantom{-}0) \\
A_3\colon & (0,\phantom{-}0,\phantom{-}1,-1,\phantom{-}0,\phantom{-}0,\phantom{-}0) \\
A_3\colon & (0,\phantom{-}1,-1,\phantom{-}0,\phantom{-}0,\phantom{-}0,\phantom{-}0)
\end{array}
$
\\\hline
$A_5 + A_1$ &
\phantom{$1.$}
$
\begin{array}{rc}
A_1\colon & (2,-1,-1,-1,-1,-1,-1) \\
A_5\colon & (0,\phantom{-}0,\phantom{-}0,\phantom{-}0,\phantom{-}0,\phantom{-}1,-1) \\
A_5\colon & (0,\phantom{-}0,\phantom{-}0,\phantom{-}0,\phantom{-}1,-1,\phantom{-}0) \\
A_5\colon & (0,\phantom{-}0,\phantom{-}0,\phantom{-}1,-1,\phantom{-}0,\phantom{-}0) \\
A_5\colon & (0,\phantom{-}0,\phantom{-}1,-1,\phantom{-}0,\phantom{-}0,\phantom{-}0) \\
A_5\colon & (0,\phantom{-}1,-1,\phantom{-}0,\phantom{-}0,\phantom{-}0,\phantom{-}0)
\end{array}
$
\\\hline
$3A_2$ &
$
\begin{array}{rcrr}
1.~A_2\colon & (1,-1,-1,-1,\phantom{-}0,\phantom{-}0,\phantom{-}0) && =: v_1 \\
1.~A_2\colon & (1,\phantom{-}0,\phantom{-}0,\phantom{-}0,-1,-1,-1) && =: v_2 \\
2.~A_2\colon & (0,\phantom{-}1,-1,\phantom{-}0,\phantom{-}0,\phantom{-}0,\phantom{-}0) && =: v_3  \\
2.~A_2\colon & (0,\phantom{-}0,\phantom{-}1,-1,\phantom{-}0,\phantom{-}0,\phantom{-}0) && =: v_4  \\
3.~A_2\colon & (0,\phantom{-}0,\phantom{-}0,\phantom{-}0,\phantom{-}1,-1,\phantom{-}0) && =: v_5  \\
3.~A_2\colon & (0,\phantom{-}0,\phantom{-}0,\phantom{-}0,\phantom{-}0,\phantom{-}1,-1) && =: v_6 
\end{array}
$
\\\hline
\end{tabular}
\end{center}\vskip-2mm
\caption{Sublattices in~$\bbZ^7$ generated by the
\mbox{$A$-,}
\mbox{$D$-,}
and
$E$-configurations}
\normalsize
\end{table}

\noindent
The~assertions now follow from mechanical~calculations.\smallskip

\noindent
In~the cases where torsion-freeness is claimed, one may easily extend the basis of the sublattice given to a basis
of~$\bbZ^7$.
For~example, consider the
types~$A_n + A_1$.
Then,~we have subsets of the lattice base consisting of
$2e_1 - e_2 - \cdots - e_7$,
$e_i - e_{i+1}$
for~$i = 3, \ldots, 6$,
$e_1$,
and~$e_7$.\smallskip

\noindent
In~the cases
$4A_1$,
$A_3 + 2A_1$,
and~$A_5 + A_1$,
the lattices may indeed be extended by the
vector~$(1,0,-1,0,-1,0,-1)$
without changing the~ranks. The~lattices obtained in this way are not further refinable
within~$\bbZ^7$.\smallskip

\noindent
In~the
case~$3A_2$,
the vector
$(v_1 + v_3 - v_4) - (v_2 + v_5 - v_6) = -3e_3 + 3e_6$
is obviously
\mbox{$3$-divisible}.
The~refined lattice has
discriminant~$3$
and is, therefore, not refinable any~further.
}
\eop
\end{theorem}

\section{Surjectivity}

\begin{coros}
\label{surj}
Let\/~$V$
be a geometrically irreducible cubic~surface
over~$\bbF_{\!q}$,
not being a~cone. If
$$\pi_V \colon \MW(V) \longrightarrow \A_0 (V^\reg)$$
is not surjective then\/
$V^\reg$~has
a nontrivial finite covering which is trivial over every\/
$\bbF_{\!q}$-rational~point.\smallskip

\noindent
{\bf Proof.}
{\em
Under~the assumption, the image of the canonical map
$V^\reg(\bbF_{\!q}) \to h_0(V^\reg)$
generates a subgroup which is not~dense. Hence,~there are
$l>1$
and a surjective, continuous homomorphism
$\alpha \colon h_0(V^\reg) \to \bbZ/l\bbZ$
sending the whole image
of~$V^\reg(\bbF_{\!q})$
to~zero.

The~same is true for the composition
$\alpha \!\circ\! \iota_{V^\reg}^\prime \colon \pi_1^t(V^\reg) \to \bbZ/l\bbZ$.
But~this simply means that the
\mbox{$l$-sheeted}
covering
of~$V^\reg$
defined
by~$\alpha \!\circ\! \iota_{V^\reg}^\prime$
has exactly
$l$
$\bbF_{\!q}$-rational
points above
every~$x \in V^\reg(\bbF_{\!q})$.
}
\eop
\end{coros}

\begin{remes}
Suppose,~$\pi_V$
were not~surjective. Then,~according to the lemma, we have a nontrivial
covering~$W$
such
that~$\#W(\bbF_{\!q}) = l \!\cdot\! \#V^\reg(\bbF_{\!q})$.
The~Weil conjectures, proven by P.\,Deligne, assure that this may be possible only for very
small~$q$.
\end{remes}

\begin{exes}
Let~$V$
be a cubic surface of
type~$4A_1$
over the finite
field~$\bbF_{\!q}$.
Then,~the canonical homomorphism
$$\pi_V \colon \MW(V) \longrightarrow \A_0 (V^\reg)$$
is surjective
for~$q > 13$.\smallskip

\noindent
{\bf Proof.}
Assume~the~contrary. Then,~according to Corollary~\ref{surj}, we have a twofold
covering~$p\colon V^\prime \to V$
ramified at the four singularities such that, over every smooth
\mbox{$\bbF_{\!q}$-rational}
point
of~$V$,
there are two
of~$V^\prime$.
Being~a cubic surface,
$V$
has at least
$q^2 - 5q + 1$~points.
Hence,~$\#V^\reg(\bbF_{\!q}) \geq q^2 -5q - 3$.

On~the other hand,
$\chi_\top(V) = 3 + 6 - 4\!\cdot\!2 = 1$,
as
$V^\reg$
is
$\Pb^2$
blown up in six points with four lines~deleted.
Therefore,~$\chi_\top (V^\prime) = 6$.
Indeed,~$V^\prime$
consists of the two sheets above 
$V^\reg$
and four points of~ramification.

We~claim~$\#V^\prime(\bbF_{\!q}) \leq q^2 + 4q + 1$.
For~this, first observe that
$V^\prime$
is simply connected as, otherwise,
$V^\reg$
had more coverings than the twofold~one.
Let~$k$
be the number of blow-ups necessary in order to
desingularize~$V^\prime$.
Then,~$\dim H^2_\et(\overline{V}, \bbQ_l) = k+4$
and one has the naive estimate
$\#\overline{V}(\bbF_{\!q}) \leq q^2 + (k+4)q + 1$.
The~claim~follows.

Consequently,
$2(q^2 -5q - 3) \leq 2 \!\cdot\! \#V^\reg(\bbF_{\!q}) \leq \#(V^\prime)^\reg (\bbF_{\!q}) \leq q^2 + 4q + 1$,
which implies
$q \leq 14$,~immediately.
\eop
\end{exes}

\section{Some observations}

\begin{lems}
Let\/~$V$
be a cubic surface over the finite
field\/~$\bbF_{\!q}$.
Suppose~that, in every equivalence class
of\/~$V^\reg(\bbF_{\!q})$,
there is a point not contained in any of the lines lying
on\/~$V$.

\begin{iii}
\item
Suppose\/~$\#\!\MW(V) = 2$.
Assume~further that not all points
of\/~$V^\reg(\bbF_{\!q})$
are contained in a~plane. Then,~for the two equivalence classes\/
$M_0, M_1$
of\/
$V^\reg(\bbF_{\!q})$,
we have the relation\/
$\#M_1 - \#M_2 = \pm q$.
\item
Suppose\/~$\#\!\MW(V) = 3$.
Then,~the three equivalence classes\/
$M_0, M_1, M_2$
of\/
$V^\reg(\bbF_{\!q})$
are of the same~size.
\end{iii}\smallskip

\noindent
{\bf Proof.}
{\em
i)
Without~restriction, assume that the classes are denoted in such a way that a line not entirely contained
in~$V$
always meets zero or two points
from~$M_1$.

Then,~fix a
point~$x \in M_1$
not contained in a line lying
on~$V$.
By~assumption, there is some
$x' \in V^\reg(\bbF_{\!q})$
outside the plane tangent
at~$x$.
The~line
$g$
connecting
$x$
and~$x'$
meets
$V$
in two distinct points
$x,y \in M_1$
and
in~$z \in M_0$.

Now,~we
intersect~$V$
with the pencil of planes
containing~$g$.
We~assert that each of the
curves~$C_t$
arising contains as many points from
$M_0$
as
from~$M_1$.
This~immediately implies the~assertion. Indeed,~equinumerosity occurs as soon as we count the points
$x$,
$y$,
and~$z$~multiply.

Let~now
$C_t$
be one of the intersection~curves. We~first observe that
$x \in C_t$~is
a smooth~point. In~fact, we do not intersect
$V$
with the tangent plane
at~$x$
since that does not
contain~$g$.
$C_t$~may
be~reducible.
However,~$x$
is, by assumption, not contained in a~line. Therefore,~for every
$p \in C_t(\bbF_{\!q})$,
there is a unique
$p' \in C_t(\bbF_{\!q})$
such that
$x$,
$p$,
and
$p'$
are~collinear.
As~$p$
and~$p'$
are in different classes, the assertion~follows.\medskip

\noindent
ii)
Here,~there are two~cases.\smallskip

\noindent
{\em First case.\/}
If~$x \in M_i$,
$y \in M_j$,
and~$z \in M_k$
are the three points of intersection of a line
with~$V$
then~$i+j+k \equiv 0 \pmod 3$.\smallskip

\noindent
We~choose a
point~$x \in M_0$
which is not contained in any of the lines
on~$V$.
Then,~for every
$p \in M_1$,
there is a unique
$p' \in M_2$
such that
$x$,
$p$,
and
$p'$
are~collinear. As~this assignment is invertible, one
has~$\#\!M_1 = \#\!M_2$.
Analogously,~a starting point
$x \in M_1$
yields the
equality~$\#\!M_2 = \#\!M_0$.\medskip

\noindent
{\em Second case.\/}
If~$x \in M_i$,
$y \in M_j$,
and~$z \in M_k$
are the three points of intersection of a line
with~$V$
then~$i+j+k \not\equiv 0 \pmod 3$.\smallskip

\noindent
We~may assume without restriction
that~$i+j+k \equiv 1 \pmod 3$.
Choose~a point
$x \in M_0$
which is not contained in any of the lines
on~$V$.
The~tangent
plane~$T_x$
contains,
besides~$x$,
only points
from~$M_1$.
Further,~there are exactly
$q+1$~of
them, as, by the assumption of this case, there is no line tangent
at~$x$
of order~three. On~the other hand,
outside~$T_x$,
the sets
$M_0$
and~$M_1$
are equinumerous since the lines
through~$x$
cause a~bijection.
Consequently,~$\#M_1 = \#M_0 + q$.

Analogously,~we obtain
$\#M_2 = \#M_1 + q$
and~$\#M_0 = \#M_2 + q$
when starting with a point
$x \in M_1$
or
$x \in M_2$,~respectively.
Thus,~the second case is~contradictory.
}
\eop
\end{lems}

\begin{defis}
Let~$V$
be a cubic surface over the finite
field~$\bbF_{\!q}$
such
that~$\#\!\MW(V) = 2$.
Then,~$V^\reg(\bbF_{\!q})$
decomposes into exactly two equivalence~classes. We~will call the equivalence class {\em negative\/} that occurs an even number of times on each~line.
\end{defis}

\begin{lemos}[{\rm Connection to the Hessian}{}]
\leavevmode\\
Let\/~$V$,
given by\/
$F(X_0, \ldots, X_3) = 0$,
be a cubic surface over the finite
field\/
$\bbF_{\!q}$
of
characteristic\/~$\neq 2$.
Suppose\/~$\#\!\MW(V) = 2$.
Then,~the following is~true.\smallskip

\noindent
If\/~$p \in V^\reg(\bbF_{\!q})$
is a negative point not lying on a line contained
in\/~$V$
then the~Hessian
$$\det\frac{\partial^2 F}{\partial X_i \partial X_j}(p)$$
is a non-square
in\/~$\bbF_{\!q}$.\smallskip

\noindent
{\bf Proof.}
{\em
Consider~the tangent
plane~$T_p$
at~$p$.
The~intersection
$C_p := V \cap T_p$
is a cubic curve with a singularity
at~$p$.
Thus,~in affine coordinates and locally
near~$p$,
the equation
of~$C_p$
is of the form
$Q(x,y) + K(x,y) = 0$
for a quadratic
form~$Q$
and a cubic
form~$K$.

By~assumption, there is no line
in~$T_p$
meeting~$p$
with
multiplicity~$3$.
This~means, in particular, that
$p \in C_p$
is a double point, not a triple~point. Further,~the two tangent directions
at~$p$
are not defined
over~$\bbF_{\!q}$.
In~other words, the binary quadratic
form~$Q$
does not represent zero
over~$\bbF_{\!q}$.
This~exactly means that minus the discriminant
of~$Q$
is a non-square
in~$\bbF_{\!q}$.
It~is a direct calculation to show that
$(-\disc Q)$
coincides, up to square factors, with the Hessian
of~$F$
at~$p$.
}
\eop
\end{lemos}

\section{Experiments}

\begin{dos}
We~let
$p$
run through the prime numbers
form~$5$
through~$101$.
For~each of the primes, we followed the classification of cubic surfaces as described in~\ref{classif} and~\ref{classifred}. For~each type, we selected ten examples by help of a random number~generator. For~those types which clearly have no moduli, we took only one~example. We~avoided the surfaces decomposing into three planes over a proper extension
of~$\bbF_{\!p}$
as, for these,
$\MW(V)$
is known to~degenerate. All~in all, we worked with 330 cubic surfaces per~prime.

For~each surface, we determined the partition
of~$V(\bbF_{\!p})$
into equivalence~classes. For~this, we run an implementation of Algorithm~\ref{alg3} in~{\tt magma}.
\end{dos}

\begin{tres}
The~partition of the points found allowed us to determine
$\MW(V)$
for every surface in the~sample.

There~is another observation, which is by far more~astonishing. In~each case, according to the theory described, we know an abelian group, 
$\MW(V)$
naturally surjects~to. It~turned out that
$\MW(V)$
was equal that group with only one~exception.

The~exception occurred
for~$p=5$.
It~was the cone over the elliptic curve given
by~$y^2 = x^3 + 2x$.
As~this elliptic curve has only two
\mbox{$\bbF_{\!5}$-rational}
points, the construction
of~$\MW(V)$
must~degenerate.
\end{tres}

\begin{remes}
This~effect clearly becomes much worse for
$p=2$
or~$3$.
This~is one of the reasons why these primes were excluded from the~experiments.
\end{remes}

\begin{suma}
\begin{iii}
\item[{\em Case I.i)\/ }]
Among the normal cubic surfaces having only double points, we always found
$\MW(V) = 0$
except for the cases
$4A_1$,
$A_3 + 2A_1$,
$A_5 + A_1$,
and~$3A_2$.
In~the first three of these cases, we have
$\MW(V) = \bbZ/2\bbZ$.

Finally,~in the
case~$3A_2$,
we established that
$\MW(V) = \bbZ/3\bbZ$
for~$p \equiv 1 \pmod 3$
and
$\Frob$
acting on the singular points by an even permutation and
for~$p \equiv 2 \pmod 3$
and
$\Frob$
acting by an odd~permutation.
Otherwise,~$\MW(V) = 0$.
\item[{\em Cases I.ii) and II.ii)\/ }]
Ignoring~the exception mentioned, for the cones,
$\MW(V)$~was
always equal to the Mordell-Weil group of the underlying~curve.
\item[{\em Case II.i)\/ }]
The~cubic ruled surfaces always
fulfilled~$\MW(V) = 0$.
\item[{\em Two components. }]
When~$V$
consisted of a non-degenerate quadric and a plane, we always found that
$\MW(V) = \bbZ$,
two points being equivalent if and only if they belonged to the same~component. When~the quadric was a cone and the plane did not meet the cusp, it turned out that
$\MW(V) = \bbZ \oplus \bbZ/2\bbZ$,
the surjection described in Example~\ref{4klassen} being~bijective.

A~cubic surface consisting of a cone and a plane through the~cusp is a cone over a reducible cubic~curve.
Here,~$\MW(V)$
was always isomorphic to the Mordell-Weil group of the~curve.
\item[{\em Three components. } ]
A~cubic surface consisting of three planes meeting in a point is the cone over a~triangle.
$\MW(V)$~was
always equal to the Mordell-Weil group of the~triangle.

Finally,~three planes meeting in a line form a cone in may~ways. Hence,~two distinct points are never equivalent to each~other. We~have
$\MW(V) \cong (K^+)^2 \oplus \bbZ^2$.
\end{iii}
\end{suma}

\begin{remes}
The~case of a cubic surface consisting of three planes with a line in common is the easiest from the theoretical point of~view. For~Algorithm~\ref{alg3}, it is, however, the most complicated~one. No~simplification occurs as no equivalent points may be~found. The~running time is dominated by steps iv) and~v), which are otherwise~negligible.
For~$p > 70$,
we excluded this case from the~experiments.
\end{remes}

\begin{remes}
On~a Quad-Core AMD Opteron Processor~2356, the total CPU~time was eight minutes
for~$p = 5$,
a little less than an hour
for~$p = 37$,
three and a half hours
for~$p = 71$,
and more than ten hours
for~$p = 101$.
\end{remes}

\appendix
\section{Algorithms}

\begin{algos}[{\rm Equivalent points}{}]
\label{alg1}
\begin{iii}
\item
Using a random number generator, choose four distinct points
$x_{11}, x_{12}, x_{21}$,
$x_{22} \in V^\reg(\bbF_{\!q})$.
\item
Determine~four points
$x_{13}, x_{23}, x_{31}, x_{32} \in V^\reg(\bbF_{\!q})$
such that the relations
$[x_{11}, x_{12}, x_{13}]$,
$[x_{21}, x_{22}, x_{23}]$,
$[x_{11}, x_{21}, x_{31}]$,
and~$[x_{12}, x_{22}, x_{32}]$
are~fulfilled. If~this turns out to be impossible as
$(x_{11}, x_{12})$,
$(x_{21}, x_{22})$,
$(x_{11}, x_{21})$,
or
$(x_{12}, x_{22})$
are lying on a line completely contained
in~$V$
then output FAIL and terminate~prematurely.
\item
Determine~points
$x_{33}$
and~$x'_{33}$
such that
$[x_{13}, x_{23}, x_{33}]$
and~$[x_{31}, x_{32}, x'_{33}]$.
If~this turns out to be impossible as
$(x_{13}, x_{23})$
or
$(x_{31}, x_{32})$
are lying on a line completely contained
in~$V$
then output FAIL and terminate~prematurely.
\item
Output
``$x_{33}$~and
$x'_{33}$
are~equivalent.''
\end{iii}
\end{algos}

\begin{algos}[{\rm A~point being equivalent to a given $x_0 \in V^\reg(\bbF_{\!q})$}{}]
\label{alg2}
\leavevmode

\begin{iii}
\item
Execute~Algorithm~\ref{alg1} in order to find two mutually equivalent points
$x_1$
and~$x_2$.
\item
Determine~a point
$x_1^\prime$
such
that~$[x_1, x_0, x_1^\prime]$.
If~this turns out to be impossible as
$(x_1, x_0)$
are lying on a line completely contained
in~$V$
then output FAIL and terminate~prematurely.
\item
Now,~determine~a point
$x_0^\prime$
such
that~$[x_1^\prime, x_2, x_0^\prime]$.
If~this turns out to be impossible as
$(x_1^\prime, x_2)$
are lying on a line completely contained
in~$V$
then output FAIL and terminate~prematurely.
\item
Output
``$x_0^\prime$~is
equivalent
to~$x_0$.''
\end{iii}
\end{algos}

\begin{algos}[{\rm Partition of the points}{}]
\label{alg3}
\leavevmode

\begin{iii}
\item
Choose~a natural
number~$N$.
\item
Decompose~$V^\reg(\bbF_{\!q})$
into a
set~${\mathfrak M} = \{N_1, \ldots, N_m\} = \{\{x_1\}, \ldots,\{x_m\}\}$
of~singletons.
\item\label{drei}
Execute~Algorithm~\ref{alg1},
$Nq^2$~times. When~two equivalent points
$x_1 \in M_k$
and~$x_2 \in M_l$
for
$k \neq l$
are found, unite
$M_k$
with~$M_l$
and reduce
$m$
by~$1$.
\item
List~the singletons still contained
in~${\mathfrak M}$,
i.e., the points that were never met in
step~\ref{drei}.
For~each element in the list obtained, execute Algorithm~\ref{alg2}
$N$~times.
When~two equivalent points
$x_1 \in M_k$
and~$x_2 \in M_l$
for
$k \neq l$
are found, unite
$M_k$
with~$M_l$
and reduce
$m$
by~$1$.
\item
If~sets of size less
than~$q$
remain
in~${\mathfrak M}$
then choose a single element from each of these~sets. For~each element in the list obtained, execute Algorithm~\ref{alg2}
$N$~times.
When~two equivalent points
$x_1 \in M_k$
and~$x_2 \in M_l$
for
$k \neq l$
are found, unite
$M_k$
with~$M_l$
and reduce
$m$
by~$1$.
\item
Output~the partition of
$V^\reg(\bbF_{\!q})$~found.
\end{iii}
\end{algos}

\begin{remss}
\begin{iii}
\item
Algorithm~\ref{alg3} finds a partition which is possibly too fine in comparison with the actual partition into equivalence~classes.
\item
In~practice, the value
$N = 7$
seems to work perfectly, for
$p = 5$
as well as for the biggest primes for which such an algorithm seems~reasonable.
\end{iii}
\end{remss}

\end{document}